November 10, 2016

# Some integrals related to the Basel problem

**Khristo N. Boyadzhiev**
Department of Mathematics and Statistics,
Ohio Northern University, Ada, OH 45810, USA
k-boyadzhiev@onu.edu

Abstract. We evaluate several arctangent and logarithmic integrals depending on a parameter. This provides a closed form summation of certain series and also integral and series representation of some classical constants.

## 1. Introduction

The famous Basel problem posed by Pietro Mengoli in 1644 and solved by Euler in 1735 asked for a closed form evaluation of the series

$$1 + \frac{1}{2^2} + \frac{1}{3^2} + \ldots$$

(see [7], [9], [10]). Euler proved that

(1) $\qquad 1 + \dfrac{1}{2^2} + \dfrac{1}{3^2} + \ldots = \dfrac{\pi^2}{6}.$

In the meantime, trying to evaluate this series, Leibniz discovered the representation

$$-\int_0^1 \frac{\log(1-t)}{t} dt = 1 + \frac{1}{2^2} + \frac{1}{3^2} + \ldots,$$

but was unable to find the numerical value of the integral (see comments in [7]). How to relate this integral to $\pi^2/6$ is discussed in [7] and [10]; it is shown that by using the complex logarithm one can solve the Basel problem. There exist, however, other integrals which can be used to quickly prove (1) without complex numbers. Possibly, the best example is the integral

(2) $$\int_0^1 \frac{\arcsin \alpha x}{\sqrt{1-x^2}} dx = \frac{1}{2}[\text{Li}_2(\alpha) - \text{Li}_2(-\alpha)]$$

for $|\alpha| \leq 1$. Here

$$\text{Li}_2(\alpha) = \sum_{n=1}^{\infty} \frac{\alpha^n}{n^2}$$

is the dilogarithm [12]. Setting $\alpha = 1$ in (2), the left hand side becomes

$$\frac{1}{2}(\arcsin x)^2 \Big|_0^1 = \frac{\pi^2}{8}$$

while the RHS is

$$\frac{1}{2}[\text{Li}_2(1) - \text{Li}_2(-1)] = \frac{3}{4}\text{Li}_2(1) = \frac{3}{4}(1 + \frac{1}{2^2} + \frac{1}{3^2} + \ldots)$$

and (1) follows immediately. This proof was recently published by Habib Bin Muzaffar in [13]. It is possibly one of the best solutions to the Basel problem. Euler also used the arcsin function in one of his proofs. Euler's approach is explained on pp. 85-86 in [15].

In this short paper we follow the idea from [13] and consider some other integrals that can be associated to the Basel problem, either solving it, or leading to similar results. In the process we evaluate a number of integrals from the tables [8] and [14].

We start with three arctangent integrals in Section 2, and then we also discuss two logarithmic integrals. Among other things, in section 2 we find the curious representation (equation (8) below)

$$\pi^3 = 192 \sum_{n=1}^{\infty} \frac{1}{n}\left(1 + \frac{1}{3} + \ldots + \frac{1}{2n-1}\right)\left(1 - \frac{1}{3} + \ldots + \frac{(-1)^{n-1}}{2n-1} - \frac{\pi}{4}\right).$$

In Section 3 we focus on several integral and series representation of some classical constants (see equations (15) and (17) below). In particular, we list two integral representations of $\zeta(3)$ and evaluate one special arctangent integral (see (16) and (19) below).

## 2. Special integrals with arctangents and logarithms

At the end of [13] it was mentioned that instead of $\arcsin(\alpha x)$ in (2) one could use $\arctan(\alpha x)$. However, as we shall see here, integrals with arctangents are not so simple.

Two natural candidates for the described method are clear

(3) $$\int_0^\infty \frac{\arctan \alpha x}{1+x^2} dx, \text{ and } \int_0^1 \frac{\arctan \alpha x}{1+x^2} dx \ .$$

When $\alpha = 1$, these integrals evaluate to a multiple of $\pi^2$ and in order to prove (1) we need to evaluate them also as a multiple of the series $\text{Li}_2(1)$. Here is how the fist one works.

**Proposition 1**. *For any $0 < \alpha < 1$,*

(4) $$2\int_0^\infty \frac{\arctan \alpha x}{1+x^2} dx = \log \alpha \, \log \frac{1-\alpha}{1+\alpha} + \text{Li}_2(\alpha) - \text{Li}_2(-\alpha) \ .$$

It easy to see (by using limits) that the RHS extends to $\alpha = 0$ and $\alpha = 1$. The function $\log \alpha \, \log \frac{1-\alpha}{1+\alpha}$ becomes zero for $\alpha = 0$ and $\alpha = 1$. With $\alpha = 1$ we compute immediately

$$(\arctan x)^2 \Big|_0^\infty = \frac{\pi^2}{4} = \frac{3}{2} \text{Li}_2(1)$$

and hence $\frac{\pi^2}{6} = \text{Li}_2(1)$.

Proof. Let $J(\alpha)$ be the LHS in (4). By differentiation (for $0 < \alpha < 1$)

$$J'(\alpha) = 2\int_0^\infty \frac{x}{(1+\alpha^2 x^2)(1+x^2)} dx = \frac{1}{1-\alpha^2} \int_0^\infty \left( \frac{2x}{1+x^2} - \frac{2\alpha^2 x}{1+\alpha^2 x^2} \right) dx$$

$$= \frac{1}{1-\alpha^2} \left[ \log \frac{1+x^2}{1+\alpha^2 x^2} \right]_0^\infty = \frac{-2\log \alpha}{1-\alpha^2} \ .$$

Thus, since $J(\alpha)$ is defined for $\alpha = 0$

$$J(\alpha) = \int_0^\alpha \frac{-2\log t}{1-t^2} dt .$$

Integrating by parts we find

$$J(\alpha) = \log \alpha \log \frac{1-\alpha}{1+\alpha} + \int_0^\alpha \left( \frac{\log(1+t)}{t} - \frac{\log(1-t)}{t} \right) dt$$

$$= \log \alpha \log \frac{1-\alpha}{1+\alpha} + \text{Li}_2(\alpha) - \text{Li}_2(-\alpha) .$$

This integral was recently evaluated in [2]. It is missing from the popular table [8], but appears in [14] as entry 2.7.4 (12). However, it appears there in a different form

$$2\int_0^\infty \frac{\arctan \alpha x}{1+x^2} dx = \frac{\pi^2}{3} - \frac{1}{2}\log^2(1+\alpha) - \text{Li}_2\left(\frac{1}{1+\alpha}\right) - \text{Li}_2(1-\alpha)$$

which is less helpful for proving (1) .

Now we look at the second integral in (3). Although it does not lead directly to the proof of (1), it provides the closed form evaluation of one interesting series.

**Proposition 2**. *For every $|\alpha| \leq 1$ we have*

(5) $$2\int_0^1 \frac{\arctan \alpha x}{1+x^2} dx = \sum_{n=0}^\infty \left(\log 2 - H_n^-\right) \frac{\alpha^{2n+1}}{2n+1} ,$$

*where*

$$H_n^- = 1 - \frac{1}{2} + \frac{1}{3} + \ldots + \frac{(-1)^{n-1}}{n}, \quad n = 1, 2, \ldots; \; H_0^- = 0,$$

*are the skew-harmonic numbers (see* [4]*).*

*In particular, for $\alpha = 1$ we find from (5)*

(6) $$\frac{\pi^2}{16} = \sum_{n=0}^\infty \left(\log 2 - H_n^-\right) \frac{1}{2n+1}.$$

Note that $H_n^-$ are the partial sums in the expansion of $\log 2$ and the series in (6) is alternating.

Proof. Using the Taylor series for $\arctan(\alpha x)$ we write

$$\int_0^1 \frac{\arctan \alpha x}{1+x^2} dx = \int_0^1 \left\{ \sum_{n=0}^{\infty} \frac{(-1)^n \alpha^{2n+1} x^{2n+1}}{2n+1} \right\} \frac{dx}{1+x^2}$$

$$= \sum_{n=0}^{\infty} \frac{(-1)^n \alpha^{2n+1}}{2n+1} \left\{ \int_0^1 \frac{x^{2n+1}}{1+x^2} dx \right\} .$$

Entry 3.241(1) in [8] says that

$$\int_0^1 \frac{x^{2n+1}}{1+x^2} dx = \frac{1}{2} \beta(n+1) ,$$

where $\beta(x)$ is the incomplete beta function (see 8.370, pp. 947-948 in [8] and also [5]). According to equation 8.375 (2) in [8] we have

$$\beta(n+1) = (-1)^n \left( \log 2 - H_n^- \right) .$$

Putting all these pieces together we arrive at (5). The proof is completed.

It would be interesting to see what happens when in Proposition 2 we replace $\arctan \alpha x$ by $(\arctan \alpha x)^2$, which leads to $\pi^3$ on the LHS.

With the notation

$$h_n = 1 + \frac{1}{3} + \frac{1}{5} + ... + \frac{1}{2n-1} \quad (n=1,2...), \; h_0 = 0$$

we have

**Proposition 3.** *For every $|\alpha| \leq 1$,*

(7) $$\int_0^1 \frac{(\arctan \alpha x)^2}{1+x^2} dx = \sum_{n=1}^{\infty} \frac{h_n}{n} \left( 1 - \frac{1}{3} + ... + \frac{(-1)^{n-1}}{2n-1} - \frac{\pi}{4} \right) \alpha^{2n} .$$

*In particular, when $\alpha = 1$, we have the curious representation*

(8) $$\pi^3 = 192 \sum_{n=1}^{\infty} \frac{h_n}{n} \left( 1 - \frac{1}{3} + ... + \frac{(-1)^{n-1}}{2n-1} - \frac{\pi}{4} \right) .$$

Proof. Using the Taylor series for $|x| \leq 1$,

$$\arctan x = x - \frac{x^3}{3} + \frac{x^5}{5} + ... + \frac{(-1)^{n-1} x^{2n-1}}{2n-1} + ... ,$$

it easy to compute the expansion

$$(\arctan \alpha x)^2 = \sum_{n=1}^{\infty} \frac{(-1)^{n-1}}{n} h_n \alpha^{2n} x^{2n} .$$

Therefore, according to entry 3.241(1) in [8]

$$\int_0^1 \frac{(\arctan \alpha x)^2}{1+x^2} dx = \sum_{n=1}^{\infty} \frac{(-1)^{n-1}}{n} h_n \alpha^{2n} \left\{ \int_0^1 \frac{x^{2n}}{1+x^2} dx \right\} .$$

$$= \sum_{n=1}^{\infty} \frac{(-1)^{n-1}}{n} h_n \alpha^{2n} \left\{ \frac{1}{2} \beta\left(n + \frac{1}{2}\right) \right\} ,$$

where, as above, $\beta(x)$ is the incomplete beta function. According to the representation

$$\beta(z) = \sum_{k=0}^{\infty} \frac{(-1)^k}{z+k}$$

(see 8.372 (1) in [8]) we compute

$$\frac{1}{2} \beta\left(n + \frac{1}{2}\right) = \frac{1}{2n+1} - \frac{1}{2n+3} + \frac{1}{2n+5} + \ldots = (-1)^n \left( \frac{\pi}{4} - 1 + \frac{1}{3} - \frac{1}{5} + \ldots + \frac{(-1)^n}{2n-1} \right) ,$$

and this finishes the proof.

Next we present two logarithmic integrals which can naturally be associated to the arctangent integrals in (3). Using the same technique, differentiation on a parameter, they can be evaluated to somewhat similar outcomes.

**Proposition 4**. *For every* $-1 < \alpha \leq 1$ ,

(9) $$\int_0^{\infty} \frac{\log(1+\alpha x)}{x(1+x)} dx = \log \alpha \log(1-\alpha) + \mathrm{Li}_2(\alpha) .$$

When $\alpha = 1$ this is the integral

$$\int_0^{\infty} \frac{\log(1+x)}{x(1+x)} dx = \frac{\pi^2}{6} ,$$

which is equivalent to entry 4.295.18 in [8] and is also a particular case of entry 2.6.10.52 in [14].

Proof. Setting $h(\alpha)$ to be the LHS we have

$$h'(\alpha) = \int_0^\infty \frac{1}{(1+\alpha x)(1+x)} dx = \left[ \frac{1}{1-\alpha} \log \frac{1+x}{1+\alpha x} \right]_0^\infty$$

$$= \frac{1}{1-\alpha} \log \frac{1}{\alpha} = \frac{-\log \alpha}{1-\alpha}.$$

From here, integrating by parts,

$$h(\alpha) = \int_0^\alpha \frac{-\log t}{1-t} dt = \log \alpha \log(1-\alpha) - \int_0^\alpha \frac{\log(1-t)}{t} dt$$

$$= \log \alpha \log(1-\alpha) + \text{Li}_2(\alpha).$$

Done!

**Proposition 5**. *For every* $-1 < \alpha \leq 1$,

(10) $$\int_0^1 \frac{\log(1+\alpha x)}{x(1+x)} dx = \text{Li}_2\left(\frac{1}{2}\right) - \text{Li}_2\left(\frac{1-\alpha}{2}\right).$$

When $\alpha = 1$, this is entry 4.291.12 in [8] and entry 2.6.10.8 in [14].

$$\int_0^1 \frac{\log(1+x)}{x(1+x)} dx = \text{Li}_2\left(\frac{1}{2}\right) = \frac{\pi^2}{12} - \frac{1}{2} \log^2 2.$$

Proof. Setting $g(\alpha)$ to be the LHS we have

$$g'(\alpha) = \int_0^1 \frac{1}{(1+\alpha x)(1+x)} dx = \frac{1}{1-\alpha} \left[ \log \frac{1+x}{1+\alpha x} \right]_0^1 = \frac{1}{1-\alpha} \log \frac{2}{1+\alpha}.$$

At the same time we notice that

$$\frac{d}{d\alpha} \text{Li}_2\left(\frac{1-\alpha}{2}\right) = \frac{-1}{1-\alpha} \sum_{n=1}^\infty \frac{1}{n}\left(\frac{1-\alpha}{2}\right)^n = \frac{1}{1-\alpha} \log\left(\frac{1+\alpha}{2}\right),$$

so the conclusion is

$$g'(\alpha) = -\frac{d}{d\alpha} \text{Li}_2\left(\frac{1-\alpha}{2}\right).$$

Therefore, for some constant $C$ we have

$$g(\alpha) = C - \text{Li}_2\left(\frac{1-\alpha}{2}\right).$$

With $\alpha = 0$ we compute $C = \text{Li}_2\left(\frac{1}{2}\right)$ and the proof is finished.

We can evaluate the integral in (10) also in terms of a power series in $\alpha$. Using the Taylor series for $\log(1+\alpha x)$ we compute

$$\int_0^1 \frac{\log(1+\alpha x)}{x(1+x)} dx = \int_0^1 \left\{\sum_{n=0}^\infty \frac{(-1)^n \alpha^{n+1} x^n}{n+1}\right\} \frac{dx}{1+x}$$

$$= \sum_{n=0}^\infty \frac{(-1)^n \alpha^{n+1}}{n+1} \left\{\int_0^1 \frac{x^n}{1+x} dx\right\} = \sum_{n=0}^\infty \frac{(-1)^n \alpha^{n+1}}{n+1} \beta(n+1)$$

$$= \sum_{n=0}^\infty \frac{(\log 2 - H_n^-)}{n+1} \alpha^{n+1}.$$

Comparing this to (10) we come to the following result.

**Corollary 6.** *For every* $-1 < \alpha \le 1$,

$$\text{Li}_2\left(\frac{1}{2}\right) - \text{Li}_2\left(\frac{1-\alpha}{2}\right) = \sum_{n=0}^\infty \frac{(\log 2 - H_n^-)}{n+1} \alpha^{n+1}.$$

*In Particular, with* $\alpha = 1$,

$$\sum_{n=0}^\infty \frac{(\log 2 - H_n^-)}{n+1} = \text{Li}_2\left(\frac{1}{2}\right) = \frac{\pi^2}{12} - \frac{1}{2}\log^2 2.$$

(See also [4].)

## 3. Two integral representations for $\zeta(3)$ and a special arctangent integral

For $|\beta| \le 1$, let $\text{Li}_p(\beta) = \sum_{n=1}^\infty \frac{\beta^n}{n^p}$ be the polylogarithm [12]

**Lema 7**. *For every* $|\beta| \le 1$ *and* $p \ge 0$,

(11) $$\int_0^1 \frac{1}{x} (\log x)^p \log \frac{1-\beta x}{1+\beta x} dx = (-1)^{p+1} \Gamma(p+1)\{\text{Li}_{p+2}(\beta) - \text{Li}_{p+2}(-\beta)\}$$

(12) $$\int_0^1 \frac{1}{x} (\log x)^p \log(1-\beta x) dx = (-1)^{p+1} \Gamma(p+1) \text{Li}_{p+2}(\beta) .$$

Equation (12) is entry 2.6.19.6 in [14].

Proof. With $x = e^{-t}$ the first integral becomes

$$(-1)^{p+1} \int_0^\infty t^p \{-\log(1-\beta e^{-t}) + \log(1+\beta e^{-t})\} dt$$

$$= (-1)^{p+1} \sum_{n=1}^\infty \frac{\beta^n}{n} \left\{ \int_0^\infty t^p e^{-nt} dt \right\} + (-1)^{p+1} \sum_{n=1}^\infty \frac{(-1)^{n-1} \beta^n}{n} \left\{ \int_0^\infty t^p e^{-nt} dt \right\}$$

$$= (-1)^{p+1} \Gamma(p+1) \left( \sum_{n=1}^\infty \frac{\beta^n}{n} \left\{ \frac{1}{n^{p+1}} \right\} + \sum_{n=1}^\infty \frac{(-1)^{n-1} \beta^n}{n} \left\{ \frac{1}{n^{p+1}} \right\} \right)$$

$$= (-1)^{p+1} \Gamma(p+1) \left( \text{Li}_{p+2}(\beta) - \text{Li}_{p+2}(-\beta) \right)$$

The same substitution in the second integral provides

$$(-1)^{p+1} \int_0^\infty t^p \{-\log(1-\beta e^{-t})\} dt = (-1)^{p+1} \sum_{n=1}^\infty \frac{\beta^n}{n} \left\{ \int_0^\infty t^p e^{-nt} dt \right\}$$

$$= (-1)^{p+1} \Gamma(p+1) \sum_{n=1}^\infty \frac{\beta^n}{n} \left\{ \frac{1}{n^{p+1}} \right\} = (-1)^{p+1} \Gamma(p+1) \text{Li}_{p+2}(\beta)$$

and the lemma is proved.

**Corollary 8**. *We have the representations*

(13) $$\zeta(3) = \frac{8}{7} \int_0^1 \arctan t \, \arctan \frac{1}{t} \frac{dt}{t} ,$$

(14) $$\zeta(3) = \int_0^1 \log(1+t) \log\left(1+\frac{1}{t}\right) \frac{dt}{t} ,$$

*where the first integral is equivalent to*

(15) $$\int_0^1 \frac{(\arctan t)^2}{t} dt = \frac{\pi}{2}G - \frac{7}{8}\zeta(3) .$$

The most remarkable integral in (15) brings together three important constants, $\pi$, the Catalan constant $G$, and $\zeta(3)$. This result is known; see entry 8 in the list [1] and also p. 18 in [6].

Proof. The starting point is equation (4), where in the integral we make the substitution $t = \alpha x$ to bring it to the form

$$2\alpha \int_0^\infty \frac{\arctan t}{\alpha^2 + t^2} dt = \log \alpha \, \log \frac{1-\alpha}{1+\alpha} + \text{Li}_2(\alpha) - \text{Li}_2(-\alpha) .$$

Here we divide both sides by $\alpha$ and integrate for $\alpha$ from 0 to 1,

$$2\int_0^\infty \arctan t \left\{ \int_0^1 \frac{d\alpha}{\alpha^2 + t^2} \right\} dt = \int_0^1 \frac{1}{\alpha} \log \alpha \, \log \frac{1-\alpha}{1+\alpha} \, d\alpha + \int_0^1 \frac{1}{\alpha} \{\text{Li}_2(\alpha) - \text{Li}_2(-\alpha)\} \, d\alpha .$$

Evaluating these integrals (using the above lemma for the second one) we come to the equation

$$\int_0^\infty \arctan t \, \arctan \frac{1}{t} \frac{dt}{t} = \text{Li}_3(1) - \text{Li}_3(-1) = \frac{7}{4}\zeta(3) .$$

We shall transform now this integral. First we split it this way: $\int_0^\infty = \int_0^1 + \int_1^\infty$, and then in the last one we make the substitution $x = \frac{1}{t}$ to get

$$\int_0^\infty \arctan t \, \arctan \frac{1}{t} \frac{dt}{t} = 2\int_0^1 \arctan t \, \arctan \frac{1}{t} \frac{dt}{t} .$$

This proves the first representation above, equation (13),

$$\int_0^1 \arctan t \, \arctan \frac{1}{t} \frac{dt}{t} = \frac{7}{8}\zeta(3) .$$

Next we use the identity $(t > 0)$

$$\arctan \frac{1}{t} = \frac{\pi}{2} - \arctan t ,$$

and the well-known fact that

$$\int_0^1 \frac{\arctan t}{t} dt = G,$$

to prove (15).

For the second representation (14) we use (9) in the form (with $t = \alpha x$)

$$\alpha \int_0^\infty \frac{\log(1+t)}{t(\alpha+t)} dt = \log \alpha \log(1-\alpha) + \text{Li}_2(\alpha).$$

Dividing by $\alpha$ and integrating for $\alpha$ from 0 to 1 we write

$$\int_0^\infty \frac{\log(1+t)}{t} \left\{ \int_0^1 \frac{d\alpha}{\alpha+t} \right\} dt = \int_0^1 \frac{1}{\alpha} \log \alpha \log(1-\alpha) d\alpha + \int_0^1 \frac{\text{Li}_2(\alpha)}{\alpha} d\alpha,$$

that is,

$$\int_0^\infty \frac{\log(1+t)}{t} \log\left(1+\frac{1}{t}\right) dt = \text{Li}_3(1) + \text{Li}_3(1) = 2\zeta(3).$$

In the same way as above we transform this integral to

$$\int_0^\infty \frac{\log(1+t)}{t} \log\left(1+\frac{1}{t}\right) dt = 2\int_0^1 \frac{\log(1+t)}{t} \log\left(1+\frac{1}{t}\right) dt,$$

which yields (14). The proof of the corollary is finished.

We end with an extension of equation (15) to power series.

**Proposition 9.** *For any* $|\alpha| \leq 1$,

(16) $$\int_0^1 \frac{(\arctan \alpha x)^2}{x} dx = \frac{1}{2} \sum_{n=1}^\infty (-1)^{n-1} \frac{h_n}{n^2} \alpha^{2n}.$$

*In particular, for* $\alpha = 1$,

(17) $$\sum_{n=1}^\infty (-1)^{n-1} \frac{h_n}{n^2} = \pi G - \frac{7}{4}\zeta(3).$$

(The numbers $h_n$ were defined right before Proposition 3.)

The series (17) is entry (59) in [6].

Proof. By expanding $(\arctan \alpha x)^2$ in power series as in the proof of Proposition 3,

$$\int_0^1 \frac{(\arctan \alpha x)^2}{x} dx = \sum_{n=1}^\infty \frac{(-1)^{n-1}}{n} h_n \alpha^{2n} \left\{\int_0^1 x^{2n-1} dx\right\} = \frac{1}{2}\sum_{n=1}^\infty \frac{(-1)^{n-1}}{n^2} h_n \alpha^{2n}.$$

With $\alpha = 1$ the assertion (17) follows from equation (15).

**Remark 10**. The series and the integral in (16) can be evaluated explicitly in a closed form by using a result of Ramanujan. Namely, Ramanujan proved that for $0 \leq \alpha \leq 1$,,

(18) $$\sum_{n=1}^\infty \frac{h_n}{n^2}\alpha^{2n} = \frac{1}{2}\log\alpha \log^2\frac{1-\alpha}{1+\alpha} + \left[\text{Li}_2\left(\frac{1-\alpha}{1+\alpha}\right) - \text{Li}_2\left(\frac{\alpha-1}{1+\alpha}\right)\right]\log\frac{1-\alpha}{1+\alpha}$$
$$- \text{Li}_3\left(\frac{1-\alpha}{1+\alpha}\right) + \text{Li}_3\left(\frac{\alpha-1}{1+\alpha}\right) + \frac{7}{4}\zeta(3)$$

(see [3], p. 255). We shall use the principle branch of the logarithm. The above equation can be extended by analytic continuation in the disc $|\alpha|<1$ where the RHS is defined. In particular, we can replace $\alpha$ by $i\alpha$ in order to obtain an alternating series. To simplify the RHS we use the dilogarithm identity

$$\text{Li}_2\left(\frac{1-\alpha}{1+\alpha}\right) - \text{Li}_2\left(\frac{\alpha-1}{1+\alpha}\right) = -\log\alpha\log\frac{1-\alpha}{1+\alpha} - \text{Li}_2(\alpha) + \text{Li}_2(-\alpha) + \frac{\pi^2}{4}$$

and also we use the formulas $\log(i\alpha) = \log\alpha + \frac{\pi}{2}i$ and $\log\frac{1-i\alpha}{1+i\alpha} = -2i\arctan\alpha$. The result is.

(19) $$\sum_{n=1}^\infty (-1)^{n-1}\frac{h_n}{n^2}\alpha^{2n} = 2\int_0^1 \frac{(\arctan\alpha x)^2}{x}dx = \text{Li}_3\left(\frac{1-i\alpha}{1+i\alpha}\right) - \text{Li}_3\left(\frac{i\alpha-1}{1+i\alpha}\right)$$
$$-2i\arctan(\alpha)\left[\text{Li}_2(i\alpha) - \text{Li}_2(-i\alpha) - \frac{\pi^2}{4}\right] - 2\left(\frac{\pi}{2}i + \log(\alpha)\right)(\arctan(\alpha))^2 - \frac{7}{4}\zeta(3).$$

**Remark 11**. Propositions 3 and 9 admit natural extensions when replacing $(\arctan \alpha x)^2$ by $(\arctan \alpha x)^p$, for any integer $p \geq 2$. In this case we use the expansion

(20) $$(\arctan x)^p = \sum_{n=1}^\infty A(n,p) x^n,$$

where $A(n,p) = 0$ for $n < p$ and for $n \geq p$

$$A(n,p) = \left( (-1)^{\frac{3n+p}{2}} + (-1)^{\frac{n-p}{2}} \right) \frac{p!}{2^{p+1}} \sum_{k=p}^{n} 2^k \binom{n-1}{k-1} \frac{s(k,p)}{k!}$$

or,

$$A(n,p) = \left( (-1)^{\frac{3n+p}{2}} + (-1)^{\frac{n-p}{2}} \right) \frac{p!}{n! 2^{p+1}} \sum_{k=p}^{n} 2^k L(n,k) s(k,p) ,$$

(see [11], Table 3). Here $s(k,p)$ are the Stirling numbers of the first kind and $L(n,k) = \binom{n-1}{k-1}\frac{n!}{k!}$ are he Lah numbers. Thus we have

(21) $$\int_0^1 \frac{(\arctan \alpha x)^p}{x} dx = \sum_{n=1}^{\infty} A(n,p) \alpha^n \left\{ \int_0^1 x^{n-1} dx \right\} = \sum_{n=1}^{\infty} A(n,p) \frac{\alpha^n}{n} ;$$

$$\int_0^1 \frac{(\arctan \alpha x)^p}{1+x^2} dx = \sum_{n=1}^{\infty} A(n,p) \alpha^n \left\{ \int_0^1 \frac{x^n}{1+x^2} dx \right\} ,$$

that is,

(22) $$\int_0^1 \frac{(\arctan \alpha x)^p}{1+x^2} dx = \frac{1}{2} \sum_{n=1}^{\infty} A(n,p) \beta\left(\frac{n+1}{2}\right) \alpha^n$$

(see [8], entry 3.241 (1)). When $\alpha = 1$ we find from here

(23) $$\pi^{p+1} = (p+1) 2^{2p+1} \sum_{n=1}^{\infty} A(n,p) \beta\left(\frac{n+1}{2}\right).$$